\documentclass{actams-author}
\usepackage{CJK}
\usepackage[all]{xy}
\numberwithin{equation}{section} 

\begin{document}

 \EditorNote{ $^*$ The author is supported by the National Natural Science Foundation of China (Grant No. 12001500, 12071444), the Natural Science Foundation of Shanxi Province of China (Grant No. 201901D111141)  and Scientific and Technological Innovation Programs of Higher Education Institutions in Shanxi (Grant No. 2020L0290).} 
\abovedisplayskip 6pt plus 2pt minus 2pt \belowdisplayskip 6pt
plus 2pt minus 2pt
\def\vsp{\vspace{1mm}}
\def\th#1{\vspace{1mm}\noindent{\bf #1}\quad}
\def\proof{\vspace{1mm}\indent{\bf Proof}\quad}
\def\no{\nonumber}
\newenvironment{prof}[1][Proof]{\indent\textbf{#1}\quad }
{\hfill $\Box$\vspace{0.7mm}}
\def\q{\quad} \def\qq{\qquad}
\allowdisplaybreaks[4]

\AuthorMark{}                             

\TitleMark{\uppercase{}}  

\title{The $\partial\bar{\partial}$-Lemma Under Surjective Maps         
}                 
\begin{CJK*}{GBK}{kai}
\author{\sl{Lingxu \uppercase{Meng}}}    
   { Department of Mathematics, North University of China,\\
 Taiyuan, Shanxi 030051,  P. R. China\\
    menglingxu@nuc.edu.cn\, }
\maketitle%
\end{CJK*}

\Abstract{We consider the $\partial\bar{\partial}$-lemma for complex manifolds under surjective holomorphic maps. Furthermore, using Deligne-Griffiths-Morgan-Sullivan's theorem, we prove that a product compact complex manifold satisfies the $\partial\bar{\partial}$-lemma if and only if so do all its components.}      


\MRSubClass{32Q99}      

\section{Introduction}  
In this paper, all complex manifolds are assumed to be \emph{connected}.

The \emph{$\partial\bar{\partial}$-lemma} on a compact complex manifold $X$ refers to that for every pure-type $d$-closed form on $X$,
the properties of $d$-exactness, $\partial$-exactness, $\bar{\partial}$-exactness and
$\partial\bar{\partial}$-exactness are equivalent while a compact complex manifold is called a \emph{$\partial\bar{\partial}$-manifold} if the \emph{$\partial\bar{\partial}$-lemma} holds on it. The most well-known examples are compact K\"{a}hler manifolds.
All complex manifolds in the Fujiki class  $\mathscr{C}$ are $\partial\bar{\partial}$-manifolds \cite{DGMS} and not K\"{a}hlerian in general.
There also exist many examples which are not in $\mathscr{C}$, see \cite{AK,AK2,ASTT,K}.
In particular, D. Angella, T. Suwa, N. Tardini and A. Tomassini \cite{ASTT} provided such an example which is simply-connected.  R. Friedman \cite{Fr} recently proved the $\partial\bar{\partial}$-lemma holds on general Clemens manifolds.

The real homotopy type of a $\partial\bar{\partial}$-manifold is a formal consequence of its de Rham cohomology ring \cite{DGMS}, which is a  strong topological obstruction of the $\partial\bar{\partial}$-lemma property.
D. Angella, N. Tardini, A. Tomassini \cite{ATar, ATom} gave criterions of $\partial\bar{\partial}$-manifolds by non-K\"{a}hlerness degrees.
The $\partial\bar{\partial}$-lemma property  is stable  under small deformations of the complex structure \cite{Vo,W,ATom}.
Recently, an important result is that the blow-up of a compact complex manifold  along a complex submanifold  is a $\partial\bar{\partial}$-manifold if and only if the complex manifold and the complex submanifold are both $\partial\bar{\partial}$-manifolds \cite{RYY, YY, ASTT, St1, St3},
which implies the heredity and bimeromorphic invariance of the $\partial\bar{\partial}$-lemma property are equivalent.

Recently, J. Stelzig  \cite{St1,St2,St3} deeply investigated the structure of bounded double complexes.
The $E_r$-isomorphism between double complexes is a central notion in his works.
In particular, the $E_1$-isomorphism is very useful for the study of various cohomologies and  the $\partial\bar{\partial}$-lemma property.
We consider the $E_1$-isomorphism under  surjective holomorphic  maps and apply it to  the $\partial\bar{\partial}$-lemma property as follows.

\begin{theorem}\label{Fr-surjective}
Let $f:X\rightarrow Y$ be a holomorphic map between compact complex manifolds  with $r=\textrm{dim}_{\mathbb{C}}X-\textrm{dim}_{\mathbb{C}}Y$.
Assume that there exists a closed $(r,r)$-current $T$ on $X$ such that $f_*T\neq0$.

$(1)$ If the Fr\"{o}licher spectral sequence  degenerates at $E_s$ for $X$, it degenerates at $E_s$ for $Y$.

$(2)$ $\Delta(X)\geq\Delta(Y)$.

$(3)$ Suppose that $X$ satisfies the $\partial\bar{\partial}$-lemma. Then $Y$ is a $\partial\bar{\partial}$-manifold.
\end{theorem}

Some applications of Theorem \ref{Fr-surjective} will be given in Sect. 3.
Furthermore, we study the $\partial\bar{\partial}$-lemma on product complex manifolds.
\begin{theorem}\label{product0}
Let $X$ and $Y$ be compact complex manifolds.
Then  $X\times Y$ satisfies the $\partial\bar{\partial}$-lemma, if and only if, $X$ and $Y$ do.
\end{theorem}

\section{Preliminaries}
\subsection{Real structure}
For an injective complex linear map  $\phi:U\rightarrow V$ of complex vector spaces, $U$ can be viewed as a complex subspace (i.e., $\phi(U)$) of $V$ via $\phi$.
Suppose that $U_1$, ..., $U_n$, $V$ are complex vector spaces and $\phi:\bigoplus\limits_{i=1}^nU_i\rightarrow V$ or $V\rightarrow \bigoplus\limits_{i=1}^nU_i$ is an isomorphism, where $\bigoplus\limits_{i=1}^nU_i$ is an external direct sum of $U_1$, ..., $U_n$.
Then $V=\bigoplus\limits_{i=1}^n\phi(U_i)$ or $V=\bigoplus\limits_{i=1}^n\phi^{-1}(U_i)$ is an internal direct sum of subspaces.
We simply write them as $V=\bigoplus\limits_{i=1}^nU_i$, which is said to be \emph{the  direct sum of the subspaces $U_1$, ..., $U_n$ via $\phi$}.

If there exsists a real vector space $K$ such that $U=K\otimes_{\mathbb{R}}\mathbb{C}$, we say that $U$ is \emph{a complex vector space with real structure}.
A complex linear map $\phi:U\rightarrow V$ of complex vector spaces  is said to be \emph{real}, if it is the complexification of a real linear map of  real vector spaces, namely, $\phi=\psi\otimes \textrm{id}:U=K\otimes_{\mathbb{R}}\mathbb{C}\rightarrow V=L\otimes_{\mathbb{R}}\mathbb{C}$ for a (real) linear map $\psi:K\rightarrow L$ of real vector spaces.
For example, $H^k(X,\mathbb{C})=H^k(X,\mathbb{R})\otimes_{\mathbb{R}}\mathbb{C}$ is a complex vector space with real structure for a smooth manifold $X$  and the pull back $f^*:H^k(Y,\mathbb{C})\rightarrow H^k(X,\mathbb{C})$ is real for a smooth map $f:X\rightarrow Y$ of smooth manifolds, since it is the complexification of $f^*:H^k(Y,\mathbb{R})\rightarrow H^k(X,\mathbb{R})$.
Assume that $\phi:U\rightarrow V$ is injective and real.
We view $U$ as a complex subspace of $V$, while the complex conjugation $\overline{W}$ of the complex subspace $W\subseteq U$ can be viewed as  the one of $\phi(W)$ in $V$ (i.e., $\overline{\phi(W)}$) via $\phi$,  since $\overline{\phi(W)}=\phi(\overline{W})$.
Suppose that $U$ and $V$ are complex vector spaces with real structures.
Then $U\otimes_{\mathbb{C}}V$ has a natural real structure, under which the complex conjugation $\overline{\alpha\otimes\beta}$ of  $\alpha\otimes\beta$ in $U\otimes_{\mathbb{C}}V$ is just $\bar{\alpha}\otimes\bar{\beta}$ for $\alpha\in U$ and $\beta\in V$.

\subsection{Double complex and  related structures}
In this subsection, we refer to \cite[Section 2]{St3} for more details.

Let $(K^{\bullet,\bullet},\partial_1,\partial_2)$ be a double complex  with two endomorphisms $\partial_1$ , $\partial_2$ of bidegree $(1,0)$ and $(0,1)$, which satisfy that $\partial_i\circ\partial_i=0$ for $i=1$, $2$ and $\partial_1\circ\partial_2+\partial_2\circ\partial_1=0$.
For convenience, we briefly write it as $K^{\bullet,\bullet}$ sometimes.
A double complex $K^{\bullet,\bullet}$ is said to be \emph{bounded}, if $K^{p,q}=0$ except for finitely many $(p,q)\in \mathbb{Z}^2$.

From now on, all the double complexes mentioned are assumed to be of vector spaces over $\mathbb{C}$ and bounded.
For a double complex $K^{\bullet,\bullet}$, let $\partial^{p,q}_1:K^{p,q}\rightarrow K^{p+1,q}$ and $\partial^{p,q}_2:K^{p,q}\rightarrow K^{p,q+1}$ be the restrictions of $\partial_1$ and $\partial_2$ respectively.
Recall the following constructions.

$\bullet$ The \emph{de Rham cohomology}
\begin{displaymath}
H_{dR}^k(K^{\bullet,\bullet}):=H^k(K^\bullet,d),
\end{displaymath}
where $(K^\bullet, d=\partial_1+\partial_2)$ is the simple complex associated to $K^{\bullet,\bullet}$.

$\bullet$ The \emph{row} and \emph{column cohomologies}
\begin{displaymath}
H^{p,q}_{\partial_1}(K^{\bullet,\bullet}):=H^{p}(K^{\bullet,q},\partial_1)\mbox{  } \mbox{  and  } \mbox{  } H^{p,q}_{\partial_2}(K^{\bullet,\bullet}):=H^{q}(K^{p,\bullet},\partial_2).
\end{displaymath}

$\bullet$ The \emph{Bott-Chern} and \emph{Aeppli cohomologies}
\begin{displaymath}
H^{p,q}_{BC}(K^{\bullet,\bullet}):=\frac{\textrm{ker}\partial^{p,q}_1\cap\textrm{ker}\partial^{p,q}_2}{\textrm{im}\partial^{p-1,q}_1\circ\partial^{p-1,q-1}_2} \mbox{  }\mbox{ and }\mbox{  }  H^{p,q}_{A}(K^{\bullet,\bullet}):=\frac{\textrm{ker}\partial^{p,q+1}_1\circ\partial^{p,q}_2}{\textrm{im}\partial^{p-1,q}_1+\textrm{im}\partial^{p,q-1}_2}.
\end{displaymath}

$\bullet$ The \emph{filtration by columns}
\begin{displaymath}
F^pK^k:=\bigoplus\limits_{\substack{r+s=k\\r\geq p}}K^{r,s}.
\end{displaymath}
This filtration on $K^\bullet$ induce the one on $H_{dR}^k(K^{\bullet,\bullet})$, i.e., the (first) \emph{Hodge filtration} on the de Rham cohomology.
We get the so called \emph{\emph{(}first\emph{)} Fr\"{o}licher spectral sequence} for $K^{\bullet,\bullet}$, denoted by $(E^{p,q}_s(K^{\bullet,\bullet}), F^pH^{k}_{dR}(K^{\bullet,\bullet}))$.
Set $b_k(K^{\bullet,\bullet}):= \textrm{dim}_{\mathbb{C}}H^{k}_{dR}(K^{\bullet,\bullet})$ and
$h_{\partial_2}^{p,q}(K^{\bullet,\bullet}):= \textrm{dim}_{\mathbb{C}}H_{\partial_2}^{p,q}(K^{\bullet,\bullet})$.
In particular, $E^{p,q}_1(K^{\bullet,\bullet})=H^{p,q}_{\partial_2}(K^{\bullet,\bullet})$
and
\begin{displaymath}
E_{\infty}^{p,k-p}(K^{\bullet,\bullet})=F^p H_{dR}^k(K^{\bullet,\bullet})/F^{p+1} H_{dR}^k(K^{\bullet,\bullet}).
\end{displaymath}
In general,
\begin{equation}\label{decrease}
\begin{aligned}
    &b_k(K^{\bullet,\bullet})=\sum\limits_{\substack{p+q=k}}\textrm{dim}_{\mathbb{C}}E_{\infty}^{p,q}(K^{\bullet,\bullet})\leq
    \ldots\leq\sum\limits_{\substack{p+q=k}}\textrm{dim}_{\mathbb{C}}E_{s}^{p,q}(K^{\bullet,\bullet})\\
\leq&\sum\limits_{\substack{p+q=k}}\textrm{dim}_{\mathbb{C}}E_{s-1}^{p,q}(K^{\bullet,\bullet})\leq\ldots   \leq\sum\limits_{\substack{p+q=k}}\textrm{dim}_{\mathbb{C}}E_{1}^{p,q}(K^{\bullet,\bullet})=\sum\limits_{\substack{p+q=k}}h_{\partial_2}^{p,q}(K^{\bullet,\bullet}).
\end{aligned}
\end{equation}
The Fr\"{o}licher spectral sequence degenerates at $E_s$ (i.e., $E_s=E_{\infty}$) if and only if $b_k(K^{\bullet,\bullet})=\sum_{\substack{p+q=k}}\textrm{dim}_{\mathbb{C}}E_{s}^{p,q}(K^{\bullet,\bullet})$ for all $k$.

Suppose that $X$ is a compact complex manifold.
Denote by $A^p(X)_{\mathbb{C}}$, $A^{p,q}(X)$, $\mathcal{D}^{p,q}(X)$ and $H^k(X,\mathbb{C})$ the spaces of complex valued smooth $p$-forms, smooth $(p,q)$-forms, $(p,q)$-currents and the de Rham cohomology with complex coefficients on  $X$, respectively.
The double complexes $(A^{\bullet,\bullet}(X),\partial,\bar{\partial})$ and  $(\mathcal{D}^{\bullet,\bullet}(X),\partial,\bar{\partial})$ are bounded.
Then $b_k(X):=b_k(A^{\bullet,\bullet}(X) )$ and $h^{p,q}(X):=h_{\bar{\partial}}^{p,q}(A^{\bullet,\bullet}(X))$ are just the $k$-th Betti and $(p,q)$-th Hodge numbers respectively.
Denote by $(E_s^{p,q}(X), F^p H^k(X,\mathbb{C}))$ the Fr\"{o}licher spectral sequence for $A^{\bullet,\bullet}(X)$.
Then $E_1^{p,q}(X)=H^{p,q}_{\bar{\partial}}(X)$ and
\begin{equation}\label{rep-Filt}
F^p H^k(X,\mathbb{C})=\{[\alpha]\in H^k(X,\mathbb{C})|\alpha\in F^pA^k(X)_{\mathbb{C}}\mbox{ and }d\alpha=0\}.
\end{equation}
Denote by $\overline{F}^qH^k(X, \mathbb{C})$  the complex conjugation of the subspace $F^qH^k(X, \mathbb{C})$ in $H^k(X, \mathbb{C})$.
Set $V^{p,q}(X):=F^pH^k(X, \mathbb{C})\cap\overline{F}^qH^k(X, \mathbb{C})$ for $p+q=k$.
We say that \emph{the Hodge filtration  gives a Hodge structure of weight $k$ on $H^k(X,\mathbb{C})$}, if
\begin{equation}\label{Hodge decom}
H^k(X, \mathbb{C})=\bigoplus_{p+q=k}V^{p,q}(X),
\end{equation}
\begin{equation}\label{sym}
\overline{V^{p,q}(X)}=V^{q,p}(X), \mbox{ for any }p+q=k.
\end{equation}

P. Deligne, P. Griffiths, J. Morgan and D. Sullivan established the well-known theorem on the relationship between the $\partial\bar{\partial}$-lemma and the Fr\"{o}licher spectral sequence as follows.

\begin{theorem}[{\cite[(5.21)]{DGMS}}]\label{DGMS}
For a compact complex manifold $X$, the following statements are equivalent:

$(i)$ The $\partial\bar{\partial}$-lemma is true for $X$.

$(ii)$ $(a)$ The Fr\"{o}licher spectral sequence of $X$ degenerates at $E_1$ $($i.e., $E_1=E_{\infty}$$)$, and

\quad \mbox{ }\mbox{ }$(b)$ the Hodge filtration  gives a Hodge structure of weight $k$ on $H^k(X,\mathbb{C})$, for every $k\geq 0$.
\end{theorem}

\begin{remark}\label{Fro-rem}
In Theorem \ref{DGMS},

$(1)$ $(ii)(a)$ is equivalent to that $F^p H^k(X,\mathbb{C})/F^{p+1} H^k(X,\mathbb{C})\cong H^{p,k-p}_{\bar{\partial}}(X)$ for all $k$, $p$, and hence is equivalent to that $b_k(X)=\sum_{\substack{p+q=k}}h^{p,q}(X)$ for all $k$.

$(2)$ $(ii)(b)$ implies that $F^pH^k(X, \mathbb{C})=\bigoplus_{\substack{r+s=k\\r\geq p}}V^{r,s}(X)$.

$(3)$ For  a $\partial\bar{\partial}$-manifold $X$, $V^{p,q}(X)\cong H_{\bar{\partial}}^{p,q}(X)$.
\end{remark}

If a morphism of double complexes induces isomorphisms on all column cohomologies (i.e., induces an isomorphism at $E_1$ between Fr\"{o}licher spectral sequences),
then it also induce  isomorphisms on all de Rham cohomologies.
A morphism of double complexes  is called an \emph{$E_1$-isomorphism}, if it induces isomorphisms on all row and column cohomologies.
We say $K^{\bullet,\bullet}\simeq_1 L^{\bullet,\bullet}$, if there exists an $E_1$-isomorphism $K^{\bullet,\bullet}\rightarrow L^{\bullet,\bullet}$.
Then $\simeq_1$ is an equivalence relation by \cite[Lemma 1.24]{St1} (or \cite[Proposition 11]{St3}).
Set $h_{BC}^{p,q}(K^{\bullet,\bullet}):= \textrm{dim}_{\mathbb{C}}H^{p,q}_{BC}(K^{\bullet,\bullet})$,
$h_{A}^{p,q}(K^{\bullet,\bullet}):= \textrm{dim}_{\mathbb{C}}H^{p,q}_{A}(K^{\bullet,\bullet})$ and
\begin{displaymath}
\Delta(K^{\bullet,\bullet}):=\sum\limits_{p,q\geq 0}h^{p,q}_{BC}(K^{\bullet,\bullet})+\sum\limits_{p,q\geq 0}h^{p,q}_{A}(K^{\bullet,\bullet})-2\sum\limits_{k\geq 0}b_{k}(K^{\bullet,\bullet}).
\end{displaymath}
By the definition, $\Delta(\bullet)$ is additive for direct sums.
Assume that $K^{\bullet,\bullet}\simeq_1 L^{\bullet,\bullet}$.
Then  $b_k(K^{\bullet,\bullet})= b_k(L^{\bullet,\bullet})$.
By \cite[Corollary 13]{St3},
$h_{BC}^{p,q}(K^{\bullet,\bullet})= h_{BC}^{p,q}(L^{\bullet,\bullet})$ and
$h_{A}^{p,q}(K^{\bullet,\bullet})= h_{A}^{p,q}(L^{\bullet,\bullet})$.
Hence, $\Delta(K^{\bullet,\bullet})=\Delta(L^{\bullet,\bullet})$.

For a compact complex manifold $X$, the inclusion $(A^{\bullet,\bullet}(X),\partial,\bar{\partial})\hookrightarrow (\mathcal{D}^{\bullet,\bullet}(X),\partial,\bar{\partial})$ is an $E_1$-isomorphism, i.e.,
$A^{\bullet,\bullet}(X)\simeq_1 \mathcal{D}^{\bullet,\bullet}(X)$.
Set $\Delta(X):=\Delta(A^{\bullet,\bullet}(X))$. D. Angella, A. Tomassini \cite[Theorems A and B]{ATom} and J. Stelzig \cite[Theorem 9]{St3} proved that
\begin{theorem}\label{ATom}
For any bounded double complex $K^{\bullet,\bullet}$, $\Delta(K^{\bullet,\bullet})\geq0$.
Moreover, a compact complex manifold $X$ satisfies the $\partial\bar{\partial}$-lemma if and only if $\Delta(X)=0$.
\end{theorem}

\section{An $E_1$-isomorphism under surjective  maps}
Consider the holomorphic map $f:X\rightarrow Y$ between compact complex manifolds with  $r=\textrm{dim}_{\mathbb{C}}X-\textrm{dim}_{\mathbb{C}}Y\geq 0$  satisfying that

$(\star)$ \emph{ there exists a closed $(r,r)$-current $T$ on $X$ such that $f_*T\neq0$.}

Notice that, \emph{any map satisfying $(\star)$ is surjective}.
Indeed, if not, by the proper mapping theorem, $f(X)$ is an analytic subset in $Y$ with dimension $<\textrm{dim}_{\mathbb{C}}Y$.
For a general point $o\in f(X)$, $X_o=f^{-1}(o)$ has dimension $\textrm{dim}_{\mathbb{C}}X-\textrm{dim}_{\mathbb{C}}f(X)> r$.
The natural map $H_{BC}^{p,q}(X)\rightarrow H^{p, q}_{BC}(\mathcal{D}^{\bullet, \bullet}(X))$ induced by the inclusion is an isomorphism, so there exists a closed smooth $(r,r)$-form $\Omega$ on $X$ such that $T=\Omega+\partial\bar{\partial}S$ for some $S\in \mathcal{D}^{r-1, r-1}(X)$.
Then $f_*T=f_*\Omega=\int_{X_o}\Omega|_{X_o}=0$, which contradicts $(\star)$.

\begin{example}\label{ex}
Let $f:X\rightarrow Y$ be a surjective  holomorphic map between compact complex manifolds with $r=\textrm{dim}_{\mathbb{C}}X-\textrm{dim}_{\mathbb{C}}Y\geq 0$.
The condition $(\star)$ was satisfied in the following cases:

$(1)$ $f$ is a  K\"ahler morphism;

$(2)$ $X$ is a $r$-K\"ahler manifold;

$(3)$ $\textrm{dim}_{\mathbb{C}}X=\textrm{dim}_{\mathbb{C}}Y$, i.e., $r=0$;

$(4)$ $X$ is a K\"ahler manifold;

$(5)$ there exists an analytic subset $Z$ of $X$ with codimension $r$ such that $f(Z)=Y$.
\end{example}
Recall that, a morphism $f:X\rightarrow Y$ of complex spaces is said to be a \emph{K\"ahler morphism}, if there exists an open covering $\mathfrak{U}=\{U_\alpha\}$ of $X$ and a system of smooth functions $\{p_\alpha\in C^\infty(U_\alpha)\}$ satisfying that the restriction $p_\alpha|_{U_\alpha\cap X_y}$ is a strictly plurisubharmonic function for every $y\in Y$ and $p_\alpha-p_\beta$ is a pluriharmonic function on $U_\alpha\cap U_\beta\cap X_y$ for every $\alpha$, $\beta$, where $X_y=f^{-1}(y)$ is the fiber of $f$ over $y$ (c.f. \cite[Definition 4.1]{Fu}).
Evidently, $\{\sqrt{-1}\partial\bar{\partial}p_\alpha\}$ give a real closed $(1,1)$-form $\omega_f$ globally defined on $X$ and the restriction $\omega_f|_{X_y}$ is a K\"ahler form on $X_y$ for every $y$.
We call such an $\omega_f$ the \emph{relative K\"ahler form} of $f$.
In  $(1)$, the constant $f_*\omega_f^r=\int_{X_y}\omega_f^r|_{X_y}> 0$ for any regular value $y$ of $f$.
A complex manifold $X$ is called \emph{$p$-K\"ahlerian}, if it admits a closed transverse $(p,p)$-form $\Omega$ (c.f. \cite[Definition 1.1, 1.3]{AB}).
In such case, $\Omega|_{Z}$ is a volume form on $Z$ for any complex submanifold $Z$ of pure dimension $p$ of $X$.
In $(2)$, there exists a closed transverse $(r,r)$-form $\Omega$ on $X$ and the constant $f_*\Omega=\int_{X_y}\Omega|_{X_y}> 0$ for any regular value $y$ of $f$.
Any complex  manifold is $0$-K\"ahlerian and any K\"ahler manifold is $p$-K\"ahlerian for $p$ smaller than its dimension,
so $(3)$ and $(4)$ are special cases of $(2)$.
In $(5)$, denote by $T_Z$ the current defined by the integral on $Z$.
Then $T_Z$ is a closed $(r,r)$-current on $X$ and $f_*T_Z=d> 0$, where $d$ is the general multiplicity of  $Z$ over $Y$.

We generalize \cite[Lemma 1.38]{St1} (or \cite[Lemma 7]{St2}) as follows.
\begin{lemma}\label{$E_1$}
Let $f:X\rightarrow Y$ be a holomorphic map between compact complex manifolds  with $r=\textrm{dim}_{\mathbb{C}}X-\textrm{dim}_{\mathbb{C}}Y$. Assume that there exists a closed $(r,r)$-current $T$ on $X$ such that $f_*T\neq0$. Then
\begin{equation}\label{E_1-surj}
(f_*\circ(T\wedge),\mbox{ }\emph{pr}):A^{\bullet,\bullet}(X)\rightarrow \mathcal{D}^{\bullet,\bullet}(Y)\oplus A^{\bullet,\bullet}(X)/f^*A^{\bullet,\bullet}(Y)
\end{equation}
is an $E_1$-isomorphism. Moreover,
\begin{displaymath}
A^{\bullet,\bullet}(X)\simeq_1 A^{\bullet,\bullet}(Y)\oplus A^{\bullet,\bullet}(X)/f^*A^{\bullet,\bullet}(Y).
\end{displaymath}
\end{lemma}
\begin{prof}
Fix an integer $p$.
Denote by $A^{p,\bullet}(X)$, $A^{\bullet,p}(X)$, $\mathcal{D}^{p,\bullet}(X)$ and $\mathcal{D}^{\bullet,p}(X)$
the complex $(A^{p,\bullet}(X),\bar{\partial})$, $(A^{\bullet,p}(X),\partial)$, $(\mathcal{D}^{p,\bullet}(X),\bar{\partial})$ and $(\mathcal{D}^{\bullet,p}(X),\partial)$, respectively.
By the assumption, $c=f_*T$ is a nonzero constant.
By the projection formula,
\begin{equation}\label{proj-form}
f_*(T\wedge f^*\alpha)=c\cdot\alpha
\end{equation}
for every $\alpha\in A^{\bullet,\bullet}(Y)$, which implies that $f^*$ is injective.
Hence, there is a short exact sequence
\begin{displaymath}
\xymatrix{
 0\ar[r] & A^{\bullet,\bullet}(Y)  \ar[r]^{f^*} & A^{\bullet,\bullet}(X)\ar[r]^{\textrm{pr}\mbox{ }\qquad}&A^{\bullet,\bullet}(X)/f^*A^{\bullet,\bullet}(Y)\ar[r] &0 }
\end{displaymath}
of double complexes, which induces a long exact sequence
\begin{equation}\label{long}
\small{
\xymatrix{
\cdots H^{q}(A^{p,\bullet}(Y))\ar[r]^{\quad f^*}&H^{q}(A^{p,\bullet}(X))\ar[r]^{\textrm{pr}\mbox{ }\qquad}&H^{q}(A^{p,\bullet}(X)/f^*A^{p,\bullet}(Y))\ar[r]^{}&H^{q+1}(A^{p,\bullet}(Y))\cdots}}
\end{equation}
for any $p$. By (\ref{proj-form}), the composition
\begin{displaymath}
\small{
\xymatrix{
 H^{q}(A^{p,\bullet}(X))\ar[r]^{T\wedge\quad} & H^{q+r}(\mathcal{D}^{p+r,\bullet}(X))\ar[r]^{\quad f_*} &H^{q}(\mathcal{D}^{p,\bullet}(Y))\ar[r]^{\cong} & H^{q}(A^{p,\bullet}(Y)) }}
\end{displaymath}
is the left inverse map of $f^*$. For all $p$, $q$, we obtain split exact sequences
\begin{displaymath}
\xymatrix{
0\ar[r]^{}& H^{q}(\mathcal{D}^{p,\bullet}(Y))\cong H^{q}(A^{p,\bullet}(Y))\ar[r]^{\qquad\qquad f^*}&H^{q}(A^{p,\bullet}(X))\ar[r]^{\textrm{pr}\qquad}&H^{q}(A^{p,\bullet}(X)/f^*A^{p,\bullet}(Y))\ar[r]^{}&0}
\end{displaymath}
from (\ref{long}), which implies that (\ref{E_1-surj}) induces isomorphisms on column cohomologies.
Similarly,
\begin{displaymath}
\xymatrix{
0\ar[r]^{}& H^{q}(\mathcal{D}^{\bullet,p}(Y))\cong H^{q}(A^{\bullet,p}(Y))\ar[r]^{\qquad\qquad f^*}&H^{q}(A^{\bullet,p}(X))\ar[r]^{\textrm{pr}\mbox{ }\qquad}&H^{q}(A^{\bullet,p}(X)/f^*A^{\bullet,p}(Y))\ar[r]^{}&0}
\end{displaymath}
are split exact sequences for all $p$, $q$ and the composition
\begin{displaymath}
\small{
\xymatrix{
 H^{q}(A^{\bullet,p}(X))\ar[r]^{T\wedge\quad} &  H^{q+r}(\mathcal{D}^{\bullet,p+r}(X))\ar[r]^{\quad f_*} &H^{q}(\mathcal{D}^{\bullet,p}(Y))\ar[r]^{\cong} & H^{q}(A^{\bullet,p}(Y)) }}
\end{displaymath}
is the left inverse map of $f^*$, which implies that (\ref{E_1-surj}) induces isomorphisms on row cohomologies.
So (\ref{E_1-surj}) is an $E_1$-isomorphism. Since $A^{\bullet,\bullet}(Y)\simeq_1 \mathcal{D}^{\bullet,\bullet}(Y)$, we get the second conclusion.
\end{prof}

Now, we give a proof of Theorem \ref{Fr-surjective} as follows.
\vspace{2mm}

$\mathbf{Proof\mbox{ }of\mbox{ }Theorem}$ \ref{Fr-surjective}\quad
Assume that the Fr\"{o}licher spectral sequence  degenerates at $E_s$ for $X$. Then $b_k(X)=\sum_{p+q=k}\textrm{dim}_{\mathbb{C}} E_s^{p,q}(X)$ for all $k$.
We have
\begin{displaymath}
\begin{aligned}
b_k(X)=&b_k(Y)+b_k(A^{\bullet,\bullet}(X)/f^*A^{\bullet,\bullet}(Y))\qquad\qquad\qquad\qquad\qquad\qquad\mbox{ }(\mbox{by Lemma \ref{$E_1$}})\\
\leq&\sum_{p+q=k}\textrm{dim}_{\mathbb{C}} E_s^{p,q}(Y)+ \sum_{p+q=k}\textrm{dim}_{\mathbb{C}} E_s^{p,q}(A^{\bullet,\bullet}(X)/f^*A^{\bullet,\bullet}(Y)) \quad(\mbox{by (\ref{decrease})})\\
= &\sum_{p+q=k}\textrm{dim}_{\mathbb{C}} E_s^{p,q}(X)  \qquad\qquad\qquad\qquad\qquad\qquad\qquad\qquad\qquad(\mbox{by Lemma \ref{$E_1$}}).
\end{aligned}
\end{displaymath}
So  $b_k(Y)=\sum_{p+q=k}\textrm{dim}_{\mathbb{C}} E_s^{p,q}(Y)$ for all $k$, which implies that the Fr\"{o}licher spectral sequence degenerates at $E_s$  for $Y$.
By Lemma \ref{$E_1$}, $\Delta(X)=\Delta(Y)+\Delta(A^{\bullet,\bullet}(X)/f^*A^{\bullet,\bullet}(Y))$.
By Theorem \ref{ATom}, $(2)$ follows.
Suppose that $X$ is a $\partial\bar{\partial}$-manifold.
By $(2)$ and Theorem \ref{ATom}, we easily get $(3)$.
Here, we give a more direct proof.
Let $\alpha$ be a $(p,q)$-form on $Y$ such that $\alpha=d\beta$ for some $\beta\in A^{p+q-1}(Y)$.
Then $f^*\alpha\in A^{p,q}(X)$ satisfies that $f^*\alpha=d(f^*\beta)$. Since $X$ is a
$\partial\bar{\partial}$-manifold, $f^*\alpha=\partial\bar{\partial}\gamma$ for some $\gamma\in A^{p-1,q-1}(X)$.
For the nonzero constant $f_*\Omega$, set $c=f_*\Omega$.
By the projection formula,
\begin{displaymath}
\alpha=c^{-1}\cdot f_*(T\wedge f^*\alpha)=\partial\bar{\partial}(c^{-1}\cdot f_*(T\wedge\gamma)).
\end{displaymath}
The natural map $H_{BC}^{p,q}(Y)\rightarrow H^{p, q}_{BC}(\mathcal{D}^{\bullet, \bullet}(Y))$ induced by the inclusion is an isomorphism, so $\alpha$ is $\partial\bar{\partial}$-exact as a smooth form on $Y$.
Therefore, $Y$ is a $\partial\bar{\partial}$-manifold.
\qed

We study the $\partial\bar{\partial}$-lemma on fiber bundles.
\begin{proposition}\label{fiberbundle}
Let $\pi:E\rightarrow X$ be a holomorphic fiber bundle  over a compact complex manifold $X$, whose  fibers are  compact complex manifolds.
Assume that there exist closed  forms $t_1$, $\ldots$, $t_r$ of pure bidegrees satisfying that the restrictions  of Dolbeault  classes $[t_1]_{\bar{\partial}}$, $\ldots$, $[t_r]_{\bar{\partial}}$ to $E_x$  is a linear basis of $H_{\bar{\partial}}^{\bullet,\bullet}(E_x)= \bigoplus_{p,q\geq0}H_{\bar{\partial}}^{p,q}(E_x)$ for every $x\in X$.

$(1)$ The Fr\"{o}licher spectral sequence  degenerates at $E_s$ for $E$, if and only if, it degenerates at $E_s$ for $Y$.

$(2)$ If the $\partial\bar{\partial}$-lemma holds on $E$, then it also holds on $X$.

$(3)$ Suppose that the restrictions  of Dolbeault  classes $[\bar{t}_1]_{\bar{\partial}}$, $\ldots$, $[\bar{t}_r]_{\bar{\partial}}$ to $E_x$  is also a linear basis of $H_{\bar{\partial}}^{\bullet,\bullet}(E_x)$. Denote by $(u_1,v_1)$, \ldots, $(u_r,v_r)$ the degrees of $t_1$, $\ldots$, $t_r$  respectively.
Then
\begin{equation}\label{E_1-fiber}
\sum\limits_{i=1}^r\pi^*(\bullet)\wedge t_i:\bigoplus_{i=1}^{r}A^{\bullet,\bullet}(X)[-u_i,-v_i]\rightarrow A^{\bullet,\bullet}(E)
\end{equation}
is an $E_1$-isomorphism and $\Delta(E)=r\cdot\Delta(X)$.
Moreover, the $\partial\bar{\partial}$-lemma holds on $E$ if and only if it  holds on $X$.
\end{proposition}
\begin{prof}
Shortly write (\ref{E_1-fiber}) as $\tau:K^{\bullet,\bullet}\rightarrow A^{\bullet,\bullet}(E)$.
Between the two double complexes, $\tau$ induces isomorphisms at $E_1$-pages
\begin{equation}\label{fibbun-Dol}
\tau^{p,q}_{\bar{\partial}}=\sum_{i=1}^{r}\pi^*(\bullet)\cup [t_i]_{\bar{\partial}}:\bigoplus_{i=1}^{r}H_{\bar{\partial}}^{p-u_i,q-v_i}(X)\rightarrow H_{\bar{\partial}}^{p,q}(E),
\end{equation}
by \cite[Theorem 1.3]{M2},
hence induces isomorphisms  $\bigoplus_{i=1}^{r}E_{s}^{p-u_i,q-v_i}(X)\tilde{\rightarrow} E_{s}^{p,q}(E)$ at $E_s$-pages
and isomorphisms $\bigoplus_{i=1}^{r}H^{k-u_i-v_i}(X,\mathbb{C})\tilde{\rightarrow} H^{k}(E,\mathbb{C})$  at $H$-pages.

Suppose that the Fr\"{o}licher spectral sequence  degenerates at $E_s$ for $X$.
Then
\begin{displaymath}
b_k(E)=\sum_{i=1}^rb_{k-u_i-v_i}(X)
=\sum_{i=1}^r\sum_{p+q=k-u_i-v_i}\textrm{dim}_{\mathbb{C}} E_s^{p,q}(X)
=\sum_{p+q=k}\textrm{dim}_{\mathbb{C}} E_s^{p,q}(E),
\end{displaymath}
which implies that the Fr\"{o}licher spectral sequence  degenerates at $E_s$  for $E$. A half of $(1)$ is proved.

Denote by $u$ the complex dimension of the fiber of $E$.
Fix a point $x\in X$.
By the hypothesis, there is some $t_i$ such that $[t_i]_{\bar{\partial}}|_{E_x}$ is the generator of $H^{u,u}_{\bar{\partial}}(E_x)\cong\mathbb{C}$.
Then $u_i=v_i=u$ and $\int_{E_x}t_i|_{E_x}\neq 0$, i.e., $\pi_*t_i\neq 0$.
By Theorem \ref{Fr-surjective} (1) (3), the other half of $(1)$ and $(2)$ follows.

Consider the commutative diagram
\begin{equation}\label{com}
\xymatrix{
 \bigoplus_{i=1}^{r}H_{\partial}^{p-u_i,q-v_i}(X) \ar[d]_{} \ar[r]^{\qquad\quad\tau^{p,q}_{\partial}}& H_{\partial}^{p,q}(E)\ar[d]^{}\\
 \bigoplus_{i=1}^{r}H_{\bar{\partial}}^{q-v_i,p-u_i}(X)      \ar[r]^{\qquad\quad\bar{\tau}^{q,p}_{\bar{\partial}}}& H_{\bar{\partial}}^{q,p}(E),  }
\end{equation}
where $\tau^{p,q}_{\partial}=\sum_{i=1}^{r}\pi^*(\bullet)\cup [t_i]_{\partial}$,
$\bar{\tau}^{q,p}_{\bar{\partial}}=\sum_{i=1}^{r}\pi^*(\bullet)\cup [\bar{t}_i]_{\bar{\partial}}$
and the two vertical maps are both the complex conjugation, i.e., mapping $[\alpha]_{\partial}$ to $[\bar{\alpha}]_{\bar{\partial}}$.
The complex conjugation is a real  isomorphism and $\bar{\tau}^{q,p}_{\bar{\partial}}$ is an (complex) isomorphism by \cite[Theorem 1.3]{M2}.
By (\ref{com}), $\tau^{p,q}_{\partial}$ is a real isomorphism and then an (complex) isomorphism.
Combining the isomorphism $\tau^{p,q}_{\bar{\partial}}$ (i.e., (\ref{fibbun-Dol})), we proved that $\tau$ is an $E_1$-isomorphism.
Since $\Delta(A^{\bullet,\bullet}(X)[-u_i,-v_i])=\Delta(X)$, $\Delta(E)=r\cdot\Delta(X)$.
By Theorem \ref{ATom}, we obtain the last conclusion.
\end{prof}

\begin{remark}
In the previous work, we proved Proposition 3.3 (3) under a stronger assumption, see \cite[Theorem 5.5]{M2}.
\end{remark}

\begin{remark}
The flag bundles satisfy the assumptions in Proposition \ref{fiberbundle} (3), see \cite[Corollary 5.7]{M2}.
\end{remark}

\begin{remark}
Fix an integer $n$. The following two statements are equivalent:

$(1)$ The $\partial\bar{\partial}$-lemma property is bimeromorphic invariant in the category of  compact complex manifolds of dimension $n$.

$(2)$ The $\partial\bar{\partial}$-lemma holds on  compact complex manifolds of dimension $n$ which are meromorphic image of $\partial\bar{\partial}$-manifolds of dimension $n$.

Assume that  $(1)$ holds.
Suppose that $X$ satisfies  the $\partial\bar{\partial}$-lemma and $f:X\dashrightarrow Y$ is a surjective meromorphic  map between  compact complex manifolds of dimension $n$.
Any desingularization $Z$ of the graph  of $f$  is a $\partial\bar{\partial}$-manifold, since it is bimeromorphic to $X$.
The map from $Z$ onto $Y$ satisfies $(\star)$ (see Example \ref{ex} $(3)$), so $Y$ satisfies the $\partial\bar{\partial}$-lemma by Theorem \ref{Fr-surjective} $(3)$.
Hence $(2)$ holds. The other part obviously follows.

In general cases, the bimeromorphic invariance of the $\partial\bar{\partial}$-lemma property is still an open problem, see \cite{RYY, YY,  ASTT, St1, St3}.
S. Yang and X.-D. Yang \cite[Theorem 1.3]{YY} proved that the $\partial\bar{\partial}$-lemma property is bimeromorphic invariant in the category of threefolds.
Hence, \emph{the $\partial\bar{\partial}$-lemma holds on the threefolds which are meromorphic images of $\partial\bar{\partial}$-threefolds}.
\end{remark}

\section{The $\partial\bar{\partial}$-lemma on product complex manifolds}
\begin{lemma}\label{product1}
Supposet that $X$ and $Y$ are compact complex manifolds.  For any $s\geq 1$, the Fr\"{o}licher spectral sequence  degenerates at $E_s$ for $X\times Y$, if and only if, it degenerates at $E_s$ for both $X$ and $Y$.
\end{lemma}
\begin{prof}
As those in the proof of Proposition \ref{fiberbundle} (3), it is easy to prove that the external product $A^{\bullet,\bullet}(X)\otimes_{\mathbb{C}}A^{\bullet,\bullet}(Y)\rightarrow A^{\bullet,\bullet}(X\times Y)$ is an $E_1$-isomorphism by algebraic K\"{u}nneth theorem and  K\"{u}nneth theorem for  Dolbeault cohomology, where $A^{\bullet,\bullet}(X)\otimes_{\mathbb{C}}A^{\bullet,\bullet}(Y)$ is the double complex associated to $A^{\bullet,\bullet}(X)$ and $A^{\bullet,\bullet}(Y)$.
So, for any $s\geq 1$, $E^{\bullet,\bullet}_s(A^{\bullet,\bullet}(X)\otimes_{\mathbb{C}}A^{\bullet,\bullet}(Y))\cong E^{\bullet,\bullet}_s(X\times Y)$,
and then
\begin{equation}\label{Kun-E_r}
E^{\bullet,\bullet}_s(X)\otimes_{\mathbb{C}} E^{\bullet,\bullet}_s(Y)\cong E^{\bullet,\bullet}_s(X\times Y)
\end{equation}
by \cite[Lemma 15]{St3}.
We have
\begin{displaymath}
\begin{aligned}
b_k(X\times Y)=&\sum_{h+l=k}b_h(X)\cdot b_l(Y)\qquad\qquad\qquad\quad\mbox{ }(\mbox{by K\"{u}nneth theorem for de Rham cohomology})\\
\leq&\sum\limits_{\substack{a+b=h\\c+d=l\\h+l=k}}\textrm{dim}_{\mathbb{C}} E_s^{a,b}(X)\cdot \textrm{dim}_{\mathbb{C}} E_s^{c,d}(Y)\mbox{ }\mbox{ }(\mbox{by (\ref{decrease})})\\
=&\sum_{p+q=k}\textrm{dim}_{\mathbb{C}} E_s^{p,q}(X\times Y).    \qquad\qquad\quad(\mbox{by (\ref{Kun-E_r})})
\end{aligned}
\end{displaymath}
Then $b_k(X\times Y)=\sum_{p+q=k}\textrm{dim}_{\mathbb{C}} E_s^{p,q}(X\times Y)$ for all $k$, if and only if, $b_k(X)=\sum_{a+b=k}\textrm{dim}_{\mathbb{C}} E_s^{a,b}(X)$ and $b_k(Y)=\sum_{c+d=k}\textrm{dim}_{\mathbb{C}} E_s^{c,d}(Y)$ for all $k$, which implies the conclusion.
\end{prof}

Now, we give a proof of Theorem \ref{product0} as follows.
\vspace{2mm}

$\mathbf{Proof\mbox{ }of\mbox{ }Theorem\mbox{ }\ref{product0}}$\quad
Let $\textrm{pr}_1$ and $\textrm{pr}_2$ be the projections of $X\times Y$ onto $X$ and $Y$ respectively.

Suppose that $X\times Y$ is a $\partial\bar{\partial}$-manifold. Fix a point $y\in Y$.
Clearly, $\textrm{pr}_1$ has a   multisection $X\times \{y\}$, so it satisfies $(\star)$ (see Example \ref{ex} (5)).
By Theorem \ref{Fr-surjective} $(3)$, $X$ satisfies the $\partial\bar{\partial}$-lemma.
Similarly, the $\partial\bar{\partial}$-lemma holds on $Y$.

Suppose that $X$ and $Y$ are both  $\partial\bar{\partial}$-manifolds.
By Theorem \ref{DGMS} and Lemma \ref{product1}, it is sufficient to check

\textbf{Claim.} The Hodge filtration  gives a Hodge structure of weight $k$ on $H^k(X\times Y,\mathbb{C})$, for every $k\geq 0$.

Denote by $\phi$ both the external products $\textrm{pr}_1^*(\bullet)\wedge \textrm{pr}_2^*(\bullet):A^\bullet(X)_{\mathbb{C}}\otimes_{\mathbb{C}}A^\bullet(Y)_{\mathbb{C}}\rightarrow A^\bullet(X\times Y)_{\mathbb{C}}$ and  $\textrm{pr}_1^*(\bullet)\cup \textrm{pr}_2^*(\bullet):H^\bullet(X,\mathbb{C})\otimes_{\mathbb{C}}H^\bullet(Y,\mathbb{C})\rightarrow H^\bullet(X\times Y,\mathbb{C})$.
Via $\phi$, we have the direct sum decomposition
\begin{equation}\label{decomposition}
\begin{aligned}
H^k(X\times Y, \mathbb{C})=&\bigoplus_{h+l=k}H^h(X,\mathbb{C})\otimes_{\mathbb{C}}H^l(Y,\mathbb{C})\qquad\quad\mbox{ }(\mbox{by K\"{u}nneth theorem})\\
=&\bigoplus_{r+s+u+v=k}V^{r,s}(X)\otimes_{\mathbb{C}} V^{u,v}(Y)\qquad\mbox{ }(\mbox{by Theorem \ref{DGMS} and (\ref{Hodge decom})}).
\end{aligned}
\end{equation}
For $a+b=p$, $h+l=k$,  $F^aA^h(X)_{\mathbb{C}}\otimes_{\mathbb{C}}F^bA^l(Y)_{\mathbb{C}}$ is  mapped into $F^pA^k(X\times Y)_{\mathbb{C}}$ by $\phi$, which implies that $F^aH^h(X, \mathbb{C})\otimes_{\mathbb{C}}F^bH^l(Y, \mathbb{C})\subseteq F^pH^k(X\times Y, \mathbb{C})$ via $\phi$ by (\ref{rep-Filt}).
For convenience, $F^pH^k(X\times Y, \mathbb{C})$ will be written as $F^pH^k$ shortly. Then
\begin{equation}\label{contain0}
\begin{aligned}
F^pH^k\supseteq&\sum_{a+b=p}\sum_{h+l=k}F^aH^h(X,\mathbb{C})\otimes_{\mathbb{C}}F^bH^l(Y,\mathbb{C})\\
=&\sum_{\substack{a+b=p, r\geq a, u\geq b\\r+s+u+v=k}}V^{r,s}(X)\otimes_{\mathbb{C}} V^{u,v}(Y).\qquad\qquad(\mbox{by Remark \ref{Fro-rem} (2)})
\end{aligned}
\end{equation}
All $V^{r,s}(X)\otimes_{\mathbb{C}} V^{u,v}(Y)$  are direct summands in the decomposition (\ref{decomposition}) of  $H^k$, so the last sum in (\ref{contain0}) is a direct sum.
Hence
\begin{equation}\label{contain}
\begin{aligned}
F^pH^k\supseteq&\bigoplus_{\substack{a+b=p, r\geq a, u\geq b\\r+s+u+v=k}}V^{r,s}(X)\otimes_{\mathbb{C}} V^{u,v}(Y)\\
=&\bigoplus\limits_{\substack{r+u\geq p\\r+s+u+v=k}}V^{r,s}(X)\otimes_{\mathbb{C}} V^{u,v}(Y).
\end{aligned}
\end{equation}
For any $p$, we have  an (non canonical) isomorphism
\begin{displaymath}
\begin{aligned}
F^pH^k/F^{p+1}H^k\cong &H_{\bar{\partial}}^{p,k-p}(X\times Y)\qquad\qquad\qquad\qquad(\mbox{by Lemma \ref{product1} and Remark \ref{Fro-rem} (1)})\\
\cong &\bigoplus\limits_{\substack{r+u=p\\s+v=k-p}}H_{\bar{\partial}}^{r,s}(X)\otimes_{\mathbb{C}} H_{\bar{\partial}}^{u,v}(Y) \qquad(\mbox{by K\"{u}nneth theorem for Dolbeault cohomology})\\
\cong &\bigoplus\limits_{\substack{r+u=p\\r+s+u+v=k}}V^{r,s}(X)\otimes_{\mathbb{C}} V^{u,v}(Y) \quad\mbox{ }(\mbox{by Remark \ref{Fro-rem} (3)}).
\end{aligned}
\end{displaymath}
and then an (non canonical) isomorphism
\begin{equation}\label{nonnatrural direct sum}
F^pH^k\cong \bigoplus\limits_{l=p}^k F^lH^k/F^{l+1}H^k
\cong \bigoplus\limits_{\substack{r+u\geq p\\r+s+u+v=k}}V^{r,s}(X)\otimes_{\mathbb{C}} V^{u,v}(Y).
\end{equation}
By (\ref{contain}) and (\ref{nonnatrural direct sum}),
$F^pH^k=\bigoplus\limits_{\substack{r+u\geq p\\r+s+u+v=k}}V^{r,s}(X)\otimes_{\mathbb{C}} V^{u,v}(Y)$.
Since $\phi$ is real,
\begin{displaymath}
\begin{aligned}
\overline{F}^qH^k=&\bigoplus\limits_{\substack{r+u\geq q\\r+s+u+v=k}}\overline{V^{r,s}(X)}\otimes_{\mathbb{C}} \overline{V^{u,v}(Y)}\\
=&\bigoplus\limits_{\substack{r+u\leq k-q\\r+s+u+v=k}}V^{r,s}(X)\otimes_{\mathbb{C}} V^{u,v}(Y) \qquad\qquad(\mbox{by Theorem \ref{DGMS} and (\ref{sym})})  .
\end{aligned}
\end{displaymath}
For $p+q=k$,
\begin{equation}\label{(p,q)-component}
V^{p,q}(X\times Y)=F^pH^k\cap\overline{F}^qH^k=\bigoplus\limits_{\substack{r+u=p\\s+v=q}}V^{r,s}(X)\otimes_{\mathbb{C}} V^{u,v}(Y).
\end{equation}
From (\ref{decomposition}) and (\ref{(p,q)-component}), the sum $\sum_{{p+q=k}}V^{p,q}(X\times Y)$ of complex subspaces in $H^k(X\times Y, \mathbb{C})$ is just the  direct sum $\bigoplus_{p+q=k}V^{p,q}(X\times Y)$ and
\begin{displaymath}
\bigoplus_{p+q=k}V^{p,q}(X\times Y)
=\bigoplus_{r+s+u+v=k}V^{r,s}(X)\otimes_{\mathbb{C}} V^{u,v}(Y)
=H^k(X\times Y, \mathbb{C}).
\end{displaymath}
Clearly, $\phi$ is real, so  $\overline{V^{p,q}(X\times Y)}=V^{q,p}(X\times Y)$ by (\ref{(p,q)-component}). We proved the claim.
\qed

\begin{remark}
For the the part of ``only if", we can also give a direct proof.
Suppose that $X\times Y$ satisfies the $\partial\bar{\partial}$-lemma. Denote $m=\textrm{dim}_{\mathbb{C}}Y$ and let $\Omega_Y$ be a closed smooth $(m,m)$-form on $Y$ such that its Dolbeault class is the inverse image of $1$  under the isomorphism $\int_Y:H_{\bar{\partial}}^{m,m}(Y)\tilde{\rightarrow}\mathbb{C}$.
Set $\Omega=\textrm{pr}_2^*\Omega_Y$.
Clearly, $\Omega$ is $d$-closed.
The constant
\begin{displaymath}
\textrm{pr}_{1*}\Omega=\int_{\{x\}\times Y}(\textrm{pr}_2^*\Omega_Y)|_{\{x\}\times Y}=\int_{Y}\Omega_Y=1,
\end{displaymath}
for any $x\in X$.
By Theorem \ref{Fr-surjective} (3), $X$ is a $\partial\bar{\partial}$-manifold.
Similarly, $Y$ satisfies the $\partial\bar{\partial}$-lemma.
\end{remark}

\begin{remark}
By Proposition \ref{fiberbundle} (3), we gave another proof in the previous work \cite[Corollary 5.8]{M2}
\end{remark}

\acknowledgements{\rm The author warmly thanks Sheng Rao, Song Yang, Xiang-Dong Yang and Jonas Stelzig \cite{St4} for their useful discussions.
The author would like to express  gratitude to referees for their helpful suggestions and careful reading of the manuscript.}

\end{document}